\documentclass[a4paper, 11pt]{article}
\usepackage{amsmath}
\usepackage{graphicx}
\usepackage{amsfonts}
\usepackage{amssymb}
\def\<{\left < }
\def\>{\right >}
\def\({\left ( }
\def\){\right )}

\def\n{\nabla}

\begin{document}
\graphicspath{{figures/}} \setlength{\baselineskip}{17pt}

\vskip .5cm
\title{\bf Ellipsoids and elliptic hyperboloids  in the  Euclidean space
 ${\Bbb E}^{n+1}$ }
\author{Dong-Soo Kim
\thanks{
  was supported by Basic Science Research Program through the National Research Foundation of Korea (NRF) funded by the Ministry of Education, Science and Technology (2010-0022926). E-mail: dosokim@chonnam.ac.kr}
\\{\small   Department of Mathematics, Chonnam National University,}
\\ {\small Kwangju 500-757, Korea}
}
\date{}
 \maketitle
\begin{center}
{\bf Abstract}
\end{center}
\par
We establish some  characterizations of elliptic hyperboloids (resp., ellipsoids) in the $(n+1)$-dimensional Euclidean space
 ${\Bbb E}^{n+1}$, using the  $n$-dimensional area of
the sections cut off by hyperplanes and
the  $(n+1)$-dimensional volume of regions between parallel hyperplanes.
We also give  a few characterizations of elliptic paraboloids  in the $(n+1)$-dimensional Euclidean space
 ${\Bbb E}^{n+1}$.
\par \vskip 0.2cm
\noindent{\bf 1. Introduction}
\par \vskip 0.2cm

In what follows we will say  that a convex  hypersurface of  ${\Bbb R}^{n+1}$ is
{\it strictly convex} if the hypersurface   is  of positive normal curvatures
 with respect to the unit normal $N$ pointing to the convex side. In particular,
 the Gauss-Kronecker curvature $K$ is  positive
 with respect to the unit normal $N$.
 We will also say that  a convex function $f:{\Bbb R}^{n}\rightarrow {\Bbb R}$ is
{\it strictly convex} if the graph of $f$ is strictly convex
with respect to the upward unit normal $N$.

\vskip0.3cm
Consider a smooth function  $g:{\Bbb R}^{n+1}\rightarrow {\Bbb R}$.
We denote by $R_g$ the set of all regular values of the function $g$.
We assume that there exists an   interval $S_g\subset R_g$ such  that for every $k\in S_g$,
the level hypersurface $M_k=g^{-1}(k)$ is a smooth strictly convex
 hypersurface  in the $(n+1)$-dimensional Euclidean space
 ${\Bbb E}^{n+1}$. We also denote by $S_g$ the maximal interval in $R_g$ which satisfies the above property.

\vskip0.3cm
 If $k\in S_g$, then we may choose a  maximal interval  $I_k\subset S_g$
  so that  each $M_{k+h}$ with $k+h\in I_k$ lies in the convex side of $M_k$. Note that
   $I_k$ is of the form $(k, a)$
   with $a>k$ or $(b, k)$ with $b<k$  according as the gradient $\nabla g$ of the function
   $g$ points to the convex side of $M_k$ or not.

\vskip0.3cm
For examples, consider two functions $g_{\pm}:{\Bbb R}^{n+1}\rightarrow {\Bbb R}$ defined by
$g(x,z)=z^2\pm (a_1^2x_1^2+\cdots +a_n^2x_n^2)$ with positive constants $a_1,\cdots,a_n$.
  Then, for the function $g_-$ we have $R_{g_-}=R-\{0\}$, $S_{g_-}=(0, \infty)$ and $I_k=(k, \infty)$, $k\in S_{g_-}$.
  For $g_+$, we  get $R_{g_+}=S_{g_+}=(0, \infty)$ and  $I_k=(0, k)$ with $k\in S_{g_+}$.
\vskip0.3cm

For a fixed point $p \in M_{k}$ with $k\in S_g$ and a sufficiently small $h$ with $k+h\in I_k$,
we consider the tangent hyperplane $\Phi$ of
 $M_{k+h}$ at some point $v\in M_{k+h}$,  which is parallel to
 the tangent hyperplane $\Psi$ of $M_k$ at $p\in M_k$.
  We denote by $A_p^*(k,h), V_p^*(k,h)$ and $S_p^*(k,h)$
 the $n$-dimensional area of the section in $\Phi$ enclosed by $\Phi \cap M_k$,
 the $(n+1)$-dimensional volume of the region bounded by $M_k$ and the hyperplane $\Phi$, and
the $n$-dimensional surface area of the region of $M_k$  between the two hyperplanes $\Phi$ and $\Psi$,
respectively.

\vskip0.3cm
In [3],  the author and Y. H. Kim studied hypersurfaces in the $(n+1)$-dimensional Euclidean space
 ${\Bbb E}^{n+1}$ defined by the graph of some function $f:{\Bbb R}^{n}\rightarrow {\Bbb R}$.
 In our notations, they proved the following characterization theorem for elliptic paraboloids
 in the $(n+1)$-dimensional Euclidean space ${\Bbb E}^{n+1}$, which extends  a result in [2].

\vskip 0.50cm

 \noindent {\bf Proposition 1.}
  Suppose that  $f:{\Bbb R}^{n}\rightarrow {\Bbb R}$ is a  strictly convex function. We consider
   the function $g:{\Bbb R}^{n+1}\rightarrow {\Bbb R}$ defined by $g(x,z)=z-f(x), x=(x_1,\cdots, x_n)$.
  Then,   the following are equivalent.

 \vskip 0.3cm
 \noindent 1) For a fixed $k\in R$, $V_p^*(k,h)$    is  a nonnegative function $\phi (h)$,
 which depends only on $h$.

 \noindent 2) For a fixed $k\in R$, $A_p^*(k,h)/|\nabla g(p)|$ is  a nonnegative function $\psi (h)$,
 which depends only on $h$. Here $\nabla g$ denotes the   gradient of $g$.

\noindent 3)  The function $f(x)$ is a quadratic polynomial given by
 $f(x)=a_1^2x_1^2+\cdots + a_n^2x_n^2$  with $a_i>0, i=1,2,\cdots, n$,
and hence every  level hypersurface $M_k$ of $g$ is an  elliptic paraboloid.

\vskip 0.3cm

Note that in the above proposition, $R_g=S_g=R$ and $I_k=(k, \infty)$.

In particular, when $n=2$, in a long series of propositions,
Archimedes  proved that every level surface $M_k$ (paraboloid of rotation) of the function $g(x,y,z)=z-a^2(x^2+y^2)$
 in the $3$-dimensional Euclidean space
 ${\Bbb E}^{3}$   satisfies $V_p^*(k,h)=ch^2$ for some constant $c$
 ([5], p.66 and Appendix A and B).
\vskip 0.50cm

In this paper, we study the  family  of strictly convex level hypersurfaces  $M_k, k\in S_g$ of a function
 $g:{\Bbb R}^{n+1}\rightarrow {\Bbb R}$ which satisfies the following conditions.

\vskip 0.30cm

\noindent  \noindent $(V^*)$:  For $k\in S_g$ with $k+h\in I_k$, $V_p^*(k,h)$ with $p\in M_k$
   is  a nonnegative function $\phi_k(h)$,
 which depends only on $k$ and $h$.
\vskip 0.30cm

 \noindent $(A^*)$: For $k\in S_g$ with $k+h\in I_k$,
  $A_p^*(k,h)/|\nabla g(p)|$ with $p\in M_k$ is  a nonnegative function $\psi_k(h)$,
 which depends only on $k$ and $h$.

\vskip 0.30cm
  \noindent $(S^*)$:  For $k\in S_g$ with $k+h\in I_k$, $S_p^*(k,h)/|\nabla g(p)|$ with $p\in M_k$ is
   a nonnegative function $\eta_k(h)$, which depends only on $k$ and $h$.

\vskip 0.3cm
As a result, first of all, we establish the following characterizations of elliptic hyperboloids.
 \vskip0.3cm

 \noindent {\bf Theorem 2.}
 Let $f:{\Bbb R}^{n}\rightarrow {\Bbb R}$ be a nonnegative  strictly convex function.
 For a nonzero real number $\alpha\in R$ with $\alpha\ne 1$,
 let's denote by  $g$ the function defined by  $g(x,z)=z^\alpha-f(x)$.
 Suppose that  the level hypersurfaces $M_k(k\in S_g)$ of $g$  in the $(n+1)$-dimensional Euclidean space
 ${\Bbb E}^{n+1}$ are strictly convex.  Then the following are equivalent.

 \noindent 1) The function $g$ satisfies  Condition $(V^*)$.

\noindent 2) The function $g$  satisfies Condition $(A^*)$.

\noindent 3)  For $k\in S_g$, $K(p)|\nabla g(p)|^{n+2}=c(k)$ is constant on $M_k$, where $K(p)$ denotes
the Gauss-Kronecker curvature of $M_k$ at $p\in M_k$  with respect to the unit normal pointing to the convex side.

\noindent 4)  The function $g$ is given  by
$$g(x,z)=z^2-(a_1^2x_1^2+\cdots + a_n^2x_n^2),$$
where $a_i>0, i=1,2,\cdots, n$.
In this case,  $R_{g}=R-\{0\}$, $S_{g}=(0, \infty)$ and $I_k=(k, \infty)$, $k\in S_{g}$.
  \vskip0.3cm

Next, in the similar way to the proof of Theorem 2, we prove the following characterizations of ellipsoids.
 \vskip0.3cm

 \noindent {\bf Theorem 3.}
 Let $f:{\Bbb R}^{n}\rightarrow {\Bbb R}$ be a nonnegative  strictly convex function.
 For a nonzero real number $\alpha\in R$ with $\alpha\ne 1$, let's denote by  $g$ the function defined by  $g(x,z)=z^\alpha+f(x)$.
 Suppose that  the level hypersurfaces $M_k(k\in S_g)$ of $g$  in the $(n+1)$-dimensional Euclidean space
 ${\Bbb E}^{n+1}$ are strictly convex.  Then the following are equivalent.

 \noindent 1) The function $g$  satisfies  Condition $(V^*)$.

\noindent 2) The function $g$  satisfies Condition $(A^*)$.

\noindent 3) For $k\in S_g$,  $K(p)|\nabla g(p)|^{n+2}=c(k)$ is constant on $M_k$, where $K(p)$ denotes
the Gauss-Kronecker curvature of $M_k$ at $p\in M_k$  with respect to the unit normal pointing to the convex side.

\noindent 4)  The function $g$ is given  by
$$g(x,z)=z^2+a_1^2x_1^2+\cdots + a_n^2x_n^2,$$
where $a_i>0, i=1,2,\cdots, n$. In this case,  $R_{g}=S_{g}=(0, \infty)$ and $I_k=(0,k)$, $k\in S_{g}$.
  \vskip0.3cm

In view of the above theorems and Lemma 9 in Section 2,  it is natural  to ask  the following question.

\vskip 0.30cm

 \noindent {\bf Question 4.} Which functions $g:{\Bbb R}^{n+1}\rightarrow {\Bbb R}$ satisfy Condition $(S^*)$?

  \vskip 0.30cm

Partially, we answer Question 4 as follows.
\vskip 0.3cm

\noindent {\bf Theorem 5.} Let $f:{\Bbb R}^{n}\rightarrow {\Bbb R}$ be a nonnegative  strictly convex function.
 For a nonzero real number $\alpha\in R$ with $\alpha\ne 1$, let's denote by  $g$ the function defined by  $g(x,z)=z^\alpha- f(x)$.
 Suppose that  the level hypersurfaces $M_k(k\in S_g)$ of $g$  in the $(n+1)$-dimensional Euclidean space
 ${\Bbb E}^{n+1}$ are strictly convex.  Then, the function $g$ does not satisfy Condition $(S^*)$.

 \vskip0.3cm

In [3], using harmonic function theory, the author and Y. H. Kim proved  Theorem 5 when  $\alpha= 1$.

 \vskip 0.30cm

Finally, we generalize the characterization theorem of [3] for  elliptic paraboloids   in the $(n+1)$-dimensional Euclidean space
 ${\Bbb E}^{n+1}$ as follows.

\vskip 0.50cm

 \noindent {\bf Theorem 6.}
  Let $f:{\Bbb R}^{n}\rightarrow {\Bbb R}$ be a nonnegative strictly convex function.
 For a nonzero real number $\alpha\in R$ with $\alpha\ne 2$, let's denote by  $g$ the function defined by  $g(x,z)=z^\alpha- f(x)$.
 Suppose that  the level hypersurfaces $M_k(k\in S_g)$ of $g$  in the $(n+1)$-dimensional Euclidean space
 ${\Bbb E}^{n+1}$ are strictly convex.  Then the following are equivalent.
\vskip0.3cm

 \noindent 1)The function $g$  satisfies  Condition $(V^*)$.

\noindent 2) The function $g$  satisfies Condition $(A^*)$.

\noindent 3) For $k\in S_g$,  $K(p)|\nabla g(p)|^{n+2}=c(k)$ is constant on $M_k$.

\noindent 4) The function $g$ is given  by
$$g(x,z)=z-a_1^2x_1^2-\cdots - a_n^2x_n^2,$$
where  $a_1, \cdots, a_n$ are positive constants. In this case, $R_g=S_g=R$ and $I_k=(k, \infty)$.

  \vskip0.3cm
 Throughout this article, all objects are smooth and connected, unless otherwise mentioned.
\vskip 0.2cm \par

\par \vskip 0.2cm
\noindent {\bf 2. Preliminaries} \vskip0.2cm
\par

Suppose that $M$ is a smooth strictly convex  hypersurface  in the $(n+1)$-dimensional Euclidean space
 ${\Bbb E}^{n+1}$ with  the unit normal $N$ pointing to the convex side.
 For a fixed point $p \in M$ and  for a sufficiently small $t>0$,
 consider the  hyperplane $\Phi$ passing through the point $p+tN(p)$ which is parallel  to
 the tangent hyperplane $\Psi$ of $M$ at $p$.

 We denote by $A_p(t), V_p(t)$ and $S_p(t)$
 the $n$-dimensional area of the section in $\Phi$ enclosed by $\Phi \cap M$,
 the $(n+1)$-dimensional volume of the region bounded by the hypersurface and the hyperplane $\Phi$ and
the $n$-dimensional surface area of the region of $M$  between the two hyperplanes $\Phi$ and $\Psi$,
respectively.

Now, we may introduce a coordinate system $(x,z)=(x_1,x_2, \cdots , x_n, z)$
 of  ${\Bbb E}^{n+1}$ with the origin $p$, the tangent hyperplane of $M$ at $p$ is the hyperplane $z=0$.
 Furthermore, we may assume that $M$ is locally  the graph of a non-negative strictly convex  function $f:{\Bbb R}^n\rightarrow {\Bbb R}$.
 Hence  $N$ is the  unit normal  pointing upward.

Then, for a sufficiently small $t>0$ we have
\smallskip
   \begin{equation}\tag{2.1}
   \begin{aligned}
   A_p(t)&= \int _{f(x)<t}1dx,
    \end{aligned}
   \end{equation}
   \begin{equation}\tag{2.2}
   \begin{aligned}
   V_p(t)&=\int _{f(x)<t}\{t-f(x)\}dx
    \end{aligned}
   \end{equation}
   and
 \begin{equation}\tag{2.3}
   S_p(t) = \int _{f(x)<t}\sqrt{1+|\nabla f|^2}dx,
   \end{equation}
 where $x=(x_1,x_2, \cdots , x_n)$, $dx=dx_1dx_2 \cdots dx_n$ and $\nabla f$ denote the gradient vector of the function $f$.

Note that  we also have
 \begin{equation}\tag{2.4}
   \begin{aligned}
   V_p(t)&=\int _{f(x)<t}\{t-f(x)\}dx\\
   &=\int _{z=0}^{t}\{\int_{f(x)<z}1dx\}dz.
       \end{aligned}
   \end{equation}
Hence, together with the fundamental theorem of calculus,  (2.4) shows that
\begin{equation}\tag{2.5}
   \begin{aligned}
  V_p'(t)=\int _{f(x)<t}1dx=A_p(t).
    \end{aligned}
   \end{equation}

In order to prove our theorems, first of all, we  need the following.
\vskip 0.50cm

 \noindent {\bf Lemma 7.} Suppose that the Gauss-Kronecker curvature
 $K(p)$ of $M$ at $p$ is positive with respect to the unit normal $N$ pointing to the convex side of $M$.
Then we have the following.

\noindent 1)
 \begin{equation}\tag{2.6}
 \lim_{t \rightarrow 0}\frac{1}{(\sqrt{t})^{n}}A_p(t)= \frac{(\sqrt{2})^{n}\omega_n}{\sqrt{K(p)}},
 \end{equation}

\noindent 2)
 \begin{equation}\tag{2.7}
 \lim_{t \rightarrow 0}\frac{1}{(\sqrt{t})^{n+2}}V_p(t)= \frac{(\sqrt{2})^{n+2}\omega_n}{(n+2)\sqrt{K(p)}},
 \end{equation}

\noindent 3)
 \begin{equation}\tag{2.8}
 \lim_{t \rightarrow 0}\frac{1}{(\sqrt{t})^{n}}S_p(t)= \frac{(\sqrt{2})^{n}\omega_n}{\sqrt{K(p)}},
 \end{equation}
where $\omega_n$ denotes the volume of  the $n$-dimensional unit ball.
\smallskip

\vskip0.3cm
\noindent {\bf Proof.} For proofs of 1), 2) and 3) with n=2, see Lemma 7 of [2]. For a  proof of 3) with arbitrary $n$, see Lemma 8 of [3].
 \vskip 0.50cm

Now, we prove the following.
 \vskip0.3cm

\noindent {\bf Lemma 8.} Consider the family of  strictly convex level hypersurfaces
 $M_k=g^{-1}(k)$ of a function $g:{\Bbb R}^{n+1}\rightarrow {\Bbb R}$ which
 satisfies Condition $(V^*)$.
 Then, for each $k\in S_g$, on $M_k$ we have
 \begin{equation} \tag{2.9}
  \begin{aligned}
K(p)|\nabla g(p)|^{n+2}=c(k),
  \end{aligned}
  \end{equation}
which is independent of $p\in M_k$,  where  $K(p)$ is the Gauss-Kronecker curvature of $M_k$ at $p$
 with respect to the unit normal $N$ pointing to the convex side and $\nabla g(p)$ denotes the gradient of $g$ at $p$.

 \vskip0.3cm
\noindent {\bf Proof.} By considering $-g$ if necessary, we may assume that $I_k$ is of the form $[k, a]$ with $a>k$, that is,
$N=\nabla g/|\nabla g|$ on $M_k$.
For a fixed point $p\in M_k$ and a  small  $t>0$,
we have
$$V_p(t)=V_p^*(k, h(t))=\phi_k(h(t)),$$
where  $h=h(t)$ is a positive function with $h(0)=0$.
 By differentiating with respect to $t$,  we get
 \begin{equation} \tag{2.10}
  \begin{aligned}
  A_p(t)=V_p'(t)=\phi_k'(h)h'(t),
  \end{aligned}
  \end{equation}
where $\phi_k'(h)$ denotes the derivative of $\phi_k$ with respect to $h$.
Hence we obtain

\begin{equation} \tag{2.11}
  \begin{aligned}
\frac{1}{(\sqrt{t})^{n}}A_p(t)=\frac{\phi_k'(h)}{(\sqrt{h})^n}(\sqrt{\frac{h(t)}{t}})^{n}h'(t).
\end{aligned}
  \end{equation}

Now we claim that
\begin{equation} \tag{2.12}
  \begin{aligned}
  \lim_{t\rightarrow 0}h'(t)=|\nabla g(p)|.
 \end{aligned}
  \end{equation}
 Assuming (2.12), we also get
\begin{equation} \tag{2.13}
  \begin{aligned}
  \lim_{t\rightarrow 0}\sqrt{\frac{h(t)}{t}}=\sqrt{|\nabla g(p)|}.
\end{aligned}
  \end{equation}
  Let us  put   $\lim_{h\rightarrow 0}\phi_k'(h)/(\sqrt{h})^n=\gamma(k)$, which is independent of $p$.
Then it follows from (2.11), (2.12), (2.13) and Lemma 7 that
 \begin{equation} \tag{2.14}
  \begin{aligned}
  K(p)|\nabla g(p)|^{n+2}=\frac{2^n\omega_n^2}{\gamma(k)^2},
\end{aligned}
  \end{equation}
which is constant on the level hypersurface  $M_k$. Thus it suffices to show that (2.12) holds.

In order to prove (2.12), we consider an orthonormal basis
$E_1, \cdots , E_n,N(p)$ of  ${\Bbb E}^{n+1}$ at $p\in M_k$, where
$E_1, \cdots , E_n\in T_p(M_k)$ and
$N(p)= \nabla g(p)/|\nabla g(p)|$ is the unit normal pointing to the convex side.
 We consider a $1$-parameter family $\Phi_t$ of hyperplanes $\Phi_t(s_1, \cdots, s_n)=p+tN(p) +\Sigma_{i=1}^{n}s_iE_i$,
 which are parallel to the tangent hyperplane $\Phi_0$ of $M_k$ at $p$.
  For small  $t>0$, there exist $s_i=s_i(t), i=1, 2, \cdots, n$ with $s_i(0)=0$
  such that $\Phi_t$ is tangent to $M_{k+h(t)}$ at $s_i=s_i(t), i=1, 2, \cdots, n$.
Hence we have
\begin{equation} \tag{2.15}
  \begin{aligned}
  k+h(t)=g(\Phi_t(s_1, \cdots, s_n))=g(p+tN(p) + \Sigma_{i=1}^{n}s_i(t)E_i).
\end{aligned}
  \end{equation}

Thus, by differentiating with respect to $t$, we get
\begin{equation} \tag{2.16}
  \begin{aligned}
  \frac{dh}{dt}=\< \nabla g(\Phi_t(s_1, \cdots, s_n)), N(p)+ \Sigma_{i=1}^{n}\frac{ds_i}{dt}E_i\>.
\end{aligned}
  \end{equation}
Note that $\Phi_t(s_1, \cdots, s_n) \rightarrow p$ as $t$ tends to $0$.
Therefore,  we obtain
\begin{equation} \tag{2.17}
  \begin{aligned}
  \lim _{t\rightarrow 0}\frac{dh}{dt}=\lim _{t\rightarrow 0}\< \nabla g(p), N(p)\>=|\nabla g(p)|,
\end{aligned}
  \end{equation}
which completes the proof of (2.12).  This completes the proof. $\square$
  \vskip 0.50cm

Similarly to the proof of Lemma 8, we may obtain
 \vskip0.3cm

\noindent {\bf Lemma 9.} Consider the family of  strictly convex level hypersurfaces
 $M_k=g^{-1}(k)$ of a function $g:{\Bbb R}^{n+1}\rightarrow {\Bbb R}$ which
 satisfies  either Condition $(A^*)$ or Condition $(S^*)$.
 Then, on $M_k$ with $k\in S_g$ we have
 \begin{equation} \tag{2.18}
  \begin{aligned}
K(p)|\nabla g(p)|^{n+2}=d(k),
  \end{aligned}
  \end{equation}
 which is independent of $p\in M_k$, where  $K(p)$ is the Gauss-Kronecker curvature of $M_k$ at $p$
 with respect to the unit normal $N$ pointing to the convex side
 and $\nabla g(p)$ denotes the gradient of $g$ at $p$.
\vskip0.3cm

\noindent {\bf Proof.} As in the proof of Lemma 8, we may assume that $I_k$ is of the form $[k, a]$ with $a>k$, that is,
$N=\nabla g/|\nabla g|$ on $M_k$. For a fixed point $p\in M_k$ and a  small  $t>0$,
we have $A_p(t)=A_p^*(k,h(t))$ for some positive function $h=h(t)$ with $h(0)=0$.
Suppose that  $M_k$
satisfies  Condition $(A^*)$. Then, we have
\begin{equation} \tag{2.19}
  \begin{aligned}
A_p(t)=A^*_p(k, h(t))=\psi_k(h(t))|\nabla g(p)|.
  \end{aligned}
  \end{equation}
Hence we obtain
\begin{equation} \tag{2.20}
  \begin{aligned}
\frac{1}{(\sqrt{t})^{n}}A_p(t)=\frac{\psi_k(h)}{(\sqrt{h})^n}(\sqrt{\frac{h(t)}{t}})^{n}|\nabla g(p)|.
\end{aligned}
  \end{equation}
 We put $\lim_{h\rightarrow 0}\psi_k(h)/(\sqrt{h})^n=\beta(k)$, which is independent of $p\in M_k$.
 Then it follows from  (2.13), (2.20) and Lemma 7 that
 \begin{equation} \tag{2.21}
  \begin{aligned}
  K(p)|\nabla g(p)|^{n+2}=\frac{2^n\omega_n^2}{\beta(k)^2},
\end{aligned}
  \end{equation}
which is constant on the level hypersurface  $M_k$.

The remaining case can be treated similarly. This completes the proof. $\square$
  \vskip 0.2cm \par

\par \vskip 0.2cm
\noindent {\bf 3. Ellipsoids and elliptic hyperboloids} \vskip0.2cm
\par

In this section, first of all, we prove Theorem 2.

For a nonzero real number $\alpha$  with $\alpha\ne 1$ and a nonnegative  convex function  $f(x)$ defined on ${\Bbb R}^n $,
 we consider the  function  $g(x,z)=z^\alpha-f(x)$.
We assume  that the  level hypersurfaces  $M_k, k\in S_g$ defined by
$g(x,z)=k$ are all strictly convex, and hence
 each $M_k, k\in S_g$ has positive  Gauss-Kronecker curvature $K$  with respect to the unit normal $N$ pointing to the convex side.

 On each $M_k$,  by differentiating,  we have for a fixed point $p=(x,z)\in M_k$,
\begin{equation} \tag{3.1}
  \begin{aligned}
 &\n f=\alpha z^{\alpha-1}\n z,\\
 &|\n g(p)|^2=\alpha^2z^{2\alpha-2}+|\nabla f(x)|^2, \\
  &z_{ij}=\frac{1}{\alpha^2z^{2\alpha-1}}(\alpha z^\alpha f_{ij}-(\alpha-1)f_if_j), \quad i,j=1,2,\cdots, n,
\end{aligned}
  \end{equation}
  where $z_i$ denotes the partial derivative of $z$ with respect to $x_i, i=1, 2, \cdots, n$, and so on.
The Gauss-Kronecker curvature $K(p)$ of $M_k$ at $p$ is given by ([6])
\begin{equation} \tag{3.2}
  \begin{aligned}
 K(p)&=\frac{\det (z_{ij})}{(\sqrt{1+|\n z|^2})^{n+2}} \\
 &=\frac{\alpha^{n+2}z^{(\alpha-1)(n+2)}\det (z_{ij})}{({\sqrt{\alpha^2z^{2\alpha-2}+|\nabla f(x)|^2}})^{n+2}},
\end{aligned}
  \end{equation}
  where the second equality follows from (3.1).
   Thus, it follows from (3.1) and (3.2) that
\begin{equation} \tag{3.3}
  \begin{aligned}
 K(p)|\n g(p)|^{n+2}&= \alpha^{n+2}z^{(\alpha-1)(n+2)}\det (z_{ij}) \\
 &=\frac{1}{\alpha^{n-2}z^{\alpha n-2\alpha +2}}{\det (\alpha z^\alpha f_{ij}-(\alpha-1)f_if_j)}.
\end{aligned}
  \end{equation}

First, suppose that  the  function  $g(x,z)=z^\alpha-f(x)$ satisfies  Condition $(V^*)$ or
 Condition $(A^*)$.
Then, it  follows from Lemma 8 or   Lemma 9 that  the  function  $g$ satisfies the condition 3) of Theorem 2.
That is, there exists a constant $c(k)$ depending on $k$ such that
\begin{equation} \tag{3.4}
  \begin{aligned}
\det (\alpha z^\alpha f_{ij}-(\alpha-1)f_if_j)=\alpha^{n-2}c(k)z^{\alpha n-2\alpha+2}.
\end{aligned}
  \end{equation}
By substituting $z^\alpha=f(x)+k$ into (3.4), we see  that  $f(x)$ satisfies
\begin{equation} \tag{3.5}
  \begin{aligned}
\det (\alpha(f(x)+k)f_{ij}-(\alpha-1)f_if_j)=\alpha^{n-2}c(k)(f(x)+k)^{n-2+2/\alpha}.
\end{aligned}
  \end{equation}

We denote by $A_i, i=1,2,\cdots,n$ the $i$-th column vector of the matrix in the left hand side of (3.5).
Then we have
\begin{equation} \tag{3.6}
  \begin{aligned}
A_i=\alpha(f(x)+k)B_i-C_i,
\end{aligned}
  \end{equation}
where
\begin{equation*}\tag{3.7}
B_i=
\begin{pmatrix}
 f_{i1} \\
 f_{i2} \\
  \vdots  \\
  f_{in}
\end{pmatrix}=\n f_i,
\quad C_i=(\alpha-1)f_i
\begin{pmatrix}
 f_{1} \\
 f_{2} \\
  \vdots \\
  f_{n}
\end{pmatrix}
=(\alpha-1)f_i\n f.
\end{equation*}
Hence, it follows from the multilinear alternating property of determinant function that
\begin{equation} \tag{3.8}
  \begin{aligned}
\det (A_1, \cdots, A_n)&= \alpha^{n}(f(x)+k)^n\det(B_1, \cdots, B_n)\\
&-\alpha^{n-1}(f(x)+k)^{n-1}\{\det(C_1,B_2, \cdots, B_n)\\
&+\cdots +\det(B_1,B_2, \cdots,B_{n-1}, C_n)\}.
\end{aligned}
  \end{equation}

Since $\det(f_{ij})=\det (B_1, \cdots, B_n)$, it follows from  (3.5) and (3.8) that
\begin{equation} \tag{3.9}
  \begin{aligned}
c(k)(f(x)+k)^{(2-\alpha)/\alpha}&=\alpha^2(f(x)+k)\det(f_{ij})-\alpha\{\det(C_1,B_2, \cdots, B_n)\\
&+\cdots +\det(B_1,B_2, \cdots,B_{n-1}, C_n)\}\\
&=A(x)k+B(x),
\end{aligned}
  \end{equation}
  where we use the following notations.
\begin{equation} \tag{3.10}
  \begin{aligned}
A(x)&=\alpha^2\det(f_{ij}),\\
B(x)&=\alpha^2f(x)\det(f_{ij})-\alpha\{\det(C_1,B_2, \cdots, B_n)\\
&+\cdots +\det(B_1,B_2, \cdots,B_{n-1}, C_n)\}.
\end{aligned}
  \end{equation}
Note that the right hand side of (3.9) is a linear polynomial in $k$
with functions in $x=(x_1, \cdots, x_n)$ as coefficients. Furthermore, note that for each $k$, $c(k)$ is
positive and $f(x)+k$ is a nonconstant function in $x$.
It follows from  (3.9) that
\begin{equation} \tag{3.11}
  \begin{aligned}
c(k)^\alpha=(A(x)k+B(x))^\alpha(k+f(x))^{\alpha-2}.
\end{aligned}
  \end{equation}

Suppose that $\alpha$ is a nonzero real number with $\alpha\ne 1,2$.
Then,  by using logarithmic differentiation of (3.11) with respect to $x_i, i=1, 2, \cdots, n$, we get
\begin{equation} \tag{3.12}
  \begin{aligned}
\alpha (\n A(x)k+\n B(x))(k+f(x))+(\alpha-2)(A(x)k+B(x))\n f(x)=0,
\end{aligned}
  \end{equation}
  which is a quadratic polynomial in $k$. It follows from (3.12) and the assumption $\alpha\ne 0,1,2$
  that $\n f(x)=0$, which is a contradiction.

Thus, by assumption, we see that $\alpha=2$ is the only possible case.
 In this case, (3.9) implies that for some constants $a$ and $b$,  $c(k)=ak+b$ with $A(x)=a$ and $B(x)=b$.
 It follows from (3.10) that
\begin{equation} \tag{3.13}
  \begin{aligned}
\det(f_{ij})=\frac{a}{4},
\end{aligned}
  \end{equation}
   and
  \begin{equation} \tag{3.14}
  \begin{aligned}
\det(C_1,B_2,& \cdots, B_n)+\cdots +\det(B_1,B_2, \cdots,B_{n-1}, C_n) =\frac{1}{2}(af(x)-b),
\end{aligned}
  \end{equation}
 where
\begin{equation} \tag{3.15}
  \begin{aligned}
B_i=\n f_i, \quad C_i=f_i \n f, \quad i=1,2,\cdots, n.
\end{aligned}
  \end{equation}

Since $f(x)$ is a nonnegative strictly convex function,  (3.13) shows
that $\det(f_{ij})$ is a positive constant  on ${\Bbb R}^n $. Hence $f(x)$ is a  quadratic polynomial
given  by ([1], [4])
\begin{equation} \tag{3.16}
  \begin{aligned}
f(x_1, \cdots, x_n)=a_1^2x_1^2+\cdots +a_n^2x_n^2, \quad a_1, \cdots, a_n >0.
\end{aligned}
  \end{equation}
 Thus,  the level hypersurfaces must be the elliptic hyperboloids
 $M_k=g^{-1}(k)$, where $g(x,z)=z^2-(a_1^2x_1^2+\cdots +a_n^2x_n^2), z>0$  with $k>0$ and $a_1, \cdots, a_n >0$.

\vskip0.3cm
Conversely, consider the function $g$ given by
$g(x,z)=z^2-f(x), z>0$ with $k>0$, where  $f(x)=a_1^2x_1^2+\cdots +a_n^2x_n^2, a_1, \cdots, a_n >0$.
For the function $g$,  we have $R_{g}=R-\{0\}$, $S_{g}=(0, \infty)$ and $I_k=(k, \infty)$, $k\in S_{g}$.

For a fixed $k>0$ and a small $h>0$, consider
 the tangent hyperplane $\Psi$ of $M_{k}$ at a point
 $p\in M_{k}$.
There exists a point $v\in M_{k+h}$ such that the tangent hyperplane $\Phi$ of $M_{k+h}$ at $v$ is parallel to the hyperplane
$\Psi$. The two points $p$ and $v$ of tangency are related by
\begin{equation} \tag{3.17}
  \begin{aligned}
v=\frac{\sqrt{k+h}}{\sqrt{k}}p, \quad p=(p_1, \cdots , p_n,\sqrt{r^2+k}), r^2=a_1^2p_1^2+\cdots+a_n^2p_n^2.
\end{aligned}
  \end{equation}
Note that  $V_p^*(k,h)$ denote the $(n+1)$-dimensional volume
of the region of $M_k$ cut off  by the  hyperplane $\Phi$.

Then the linear mapping
\begin{equation}\tag{3.18}
   \begin{aligned}
   T_1(x_1,x_2,\cdots , x_n,z)=(a_1x_1,a_2x_2,\cdots , a_nx_n, z)
    \end{aligned}
   \end{equation}
transforms $M_k$ (resp., $M_{k+h}$) onto a hyperboloid of revolution $M'_k:z^2=x_1^2+x_2^2+\cdots +x_n^2+k,$
(resp., $M_{k+h}':z^2=x_1^2+x_2^2+\cdots +x_n^2+k+h $), $\Phi$ to a hyperplane $\Phi'$,
$p\in M_{k}$ and $v\in M_{k+h}$ to points of tangency $p' =(p_1', \cdots, p_n',\sqrt{(r')^2+k})\in M'_{k}$
and $v' =(\sqrt{k+h}/\sqrt{k})p'\in M'_{k+h}$, respectively,
where $(r')^2=\Sigma (p_i')^2$.
If we let $V_{p'}'^*(k,h)$ denote the volume  of the region of $M_k'$ cut off  by the  hyperplane $\Phi'$,
then we get
\begin{equation}\tag{3.19}
   \begin{aligned}
   V_{p'}'^*(k,h)=a_1\cdots a_n  V_p^*(k,h).
    \end{aligned}
   \end{equation}

Let's consider the rotation $A$ around the $z$-axis which maps the point $p'$ of tangency to
$p''=(0,\cdots,0,r', \sqrt{(r')^2+k})$.
Then the rotation $A$ takes $v'$ to $v''=(\sqrt{k+h}/\sqrt{k})p''$.
Note that the 1-parameter group $B(t)$ on the $x_nz$-plane defined by
\begin{equation*}\tag{3.20}
B(t)=
\begin{pmatrix}
 \cosh t& \sinh t \\
 \sinh t& \cosh t \\
 \end{pmatrix},
\end{equation*}
takes the upper hyperbola  $z^2=x_n^2+k, z>0$ (resp., $z^2=x_n^2+k+h, z>0$) onto itself.
Hence, there exists a parameter $t_0$ such that $B(t_0)$ maps $p''$ to $p'''=(0,\dots,0, \sqrt{k})$
(resp., $v''$ to $v'''=(0,\dots,0, \sqrt{k+h}))$.

We  consider the linear mapping $T_2=\bar B(t_0)\circ A$ of ${\Bbb R}^{n+1} $, where $\bar B(t_0)$ denotes the
extended linear mapping of  $B(t_0)$ on  ${\Bbb R}^{n+1} $  fixing $x_1\cdots x_{n-1}$-plane.
Then  the linear mapping $T_2$
takes the hyperboloid of rotation $M_k'$ (resp., $M_{k+h}'$)  onto itself, $p'$ and $v'$ to the points  of tangency
$p'''=(0,\dots,0, \sqrt{k})$ and $v'''=(0,\dots,0, \sqrt{k+h})$,
$\Phi'$ to the hyperplane
$\Phi'':z=\sqrt{k+h}$.

Due to the volume-preserving property of  $T_2$, we obtain
\begin{equation}\tag{3.21}
   \begin{aligned}
   V_{p'}'^*(k,h)=V^*(k,h),
    \end{aligned}
   \end{equation}
where $V^*(k,h)$ denotes  the volume  of the region of $M_k'$ cut off  by the  hyperplane $\Phi''$.
Together with (3.19), it follows from (3.21) that
\begin{equation}\tag{3.22}
   \begin{aligned}
 V_p^*(k,h)=\frac{\omega_n}{a_1\cdots a_n}\{\sqrt{k+h}h^{n/2}- n\int^{\sqrt{h}}_0\sqrt{r^2+k}r^{n-1}dr\},
    \end{aligned}
   \end{equation}
where $\omega_{n}$ denotes the volume of the $n$-dimensional unit ball. Hence, we see that
$V_p^*(k,h)$ is independent of the point $p\in M_k$, which is denoted by $\phi_k(h)$.
Thus  the function $g$ given by
$g(x,z)=z^2-f(x), z>0$
satisfies Condition $(V^*)$.

\vskip 0.30cm
Finally, we  show that the function $g$ given by
$g(x,z)=z^2-f(x), z>0$ with $k>0$,
where  $f(x)=a_1^2x_1^2+\cdots +a_n^2x_n^2, a_1, \cdots, a_n >0$
satisfies Condition $(A^*)$. For a fixed point $p=(p_1, \cdots, p_n, \sqrt {r^2+k})\in M_k$,
where $r^2=a_1^2p_1^2+ \cdots +a_n^2p_n^2$, and a  small  $t\in R$,
we have $V_p(t)=\phi_k(h(t))$ for some $h=h(t)$ with $h(0)=0$.
 By differentiating with respect to $t$,  from (2.5) we get
 \begin{equation} \tag{3.23}
  \begin{aligned}
  A_p(t)=V_p'(t)=\phi_k'(h)h'(t),
  \end{aligned}
  \end{equation}
where $\phi_k'(h)$ denotes the derivative of $\phi_k$ with respect to $h$.

With the aid of (3.17), it is straightforward to show that
 \begin{equation} \tag{3.24}
  \begin{aligned}
 t=\frac{2\sqrt{k}(\sqrt{k+h}-\sqrt{k})}{|\nabla g(p)|}.
  \end{aligned}
  \end{equation}
Hence we get
\begin{equation} \tag{3.25}
  \begin{aligned}
 h(t)=\frac{1}{4k}|\nabla g(p)|^2t^2 + |\nabla g(p)|t.
  \end{aligned}
  \end{equation}
It follows from (3.24), (3.25) and (3.23) that
\begin{equation} \tag{3.26}
  \begin{aligned}
  A_p^*(k,h)=A_p(t)=\frac{\sqrt{k+h}}{\sqrt{k}}\phi_k'(h)|\nabla g(p)|,
  \end{aligned}
  \end{equation}
which shows that the function $g$ given by
$g(x,z)=z^2-f(x), z>0$ satisfies Condition $(A^*)$.

 It follows from (3.9) and (3.10) with $\alpha =2$  that the function $g$ given by
$g(x,z)=z^2-f(x), z>0$ satisfies
   \begin{equation}\tag{3.27}
   \begin{aligned}
 K(p)|\nabla g(p)|^{n+2}=c(k)=2^{n+2}a_1^2a_2^2\cdots a_n^2k.
    \end{aligned}
   \end{equation}

This completes the proof of Theorem 2.
 \vskip 0.30cm

 Next, we prove  Theorem 3 as follows.

 For a nonzero real number $\alpha$  with $\alpha\ne 1$ and a nonnegative  convex function  $f(x)$ defined on ${\Bbb R}^n $,
 we consider the  function  $g(x,z)=z^\alpha+f(x)$.
We assume  that the  level hypersurfaces  $M_k, k\in S_g$ defined by
$g(x,z)=k$ are all strictly convex, and hence each $M_k, k\in S_g$ has positive
Gauss-Kronecker curvature $K$ with respect to the unit normal $N$ pointing to the convex side.

Suppose that  the  function  $g$ satisfies  Condition $(V^*)$ or
 Condition $(A^*)$.
Then, it  follows from Lemma 8 or   Lemma 9 that  $M_k$ satisfies the condition 3) of Theorem 3.
Then, changing $f(x)$ by  $-f(x)$ in the proof of Theorem 2,  (3.11) shows that $f(x)$ satisfies
\begin{equation} \tag{3.28}
  \begin{aligned}
(-1)^{\alpha n}c(k)^\alpha=(A(x)k+B(x))^\alpha(k-f(x))^{\alpha-2},
\end{aligned}
  \end{equation}
where
\begin{equation} \tag{3.29}
  \begin{aligned}
A(x)&=\alpha^2\det(f_{ij}),\\
B(x)&=-\alpha^2f(x)\det(f_{ij})+\alpha\{\det(C_1,B_2, \cdots, B_n)\\
&+\cdots +\det(B_1,B_2, \cdots,B_{n-1}, C_n)\},\\
\end{aligned}
  \end{equation}
and $B_i=\nabla f_i, C_i=(\alpha -1)f_i\nabla f, i=1,2,\cdots, n$.
Since $c(k)>0$, the logarithmic differentiation of (3.28) shows that $\alpha=2$, and hence
for some constants $a$ and $b$, we get $(-1)^nc(k)=ak+b$ with $\det (f_{ij})=a/4$.
Since $f(x)$ is a nonnegative strictly convex function,  this implies
that $\det(f_{ij})$ is a positive constant  on ${\Bbb R}^n $.
By the same argument as in the proof of Theorem 2,
we see that  $f(x)$ is a  quadratic polynomial
given  by
\begin{equation} \tag{3.30}
  \begin{aligned}
f(x_1, \cdots, x_n)=a_1^2x_1^2+\cdots +a_n^2x_n^2, \quad a_1, \cdots, a_n >0.
\end{aligned}
  \end{equation}
 Thus,  the level hypersurfaces must be the ellipsoids given by $g(x,z)=z^2+a_1^2x_1^2+\cdots +a_n^2x_n^2=k$
 with $k>0$ and $a_1, \cdots, a_n >0$.

\vskip0.3cm

Conversely, we consider the function $g$ given by $g(x,z)=z^2+a_1^2x_1^2+\cdots +a_n^2x_n^2$
 with  $a_1, \cdots, a_n >0$. For the function $g$, we have $R_{g}=S_{g}=(0, \infty)$ and $I_k=(0,k)$, $k\in S_{g}$.

For a fixed $k>0$ and a point $p\in M_k$, consider
 the tangent hyperplane $\Psi$ of $M_{k}$ at
 $p$.
For a sufficiently small $h<0$ with $k+h>0$,
there exists a point $v\in M_{k+h}$ such that the tangent hyperplane $\Phi$ of $M_{k+h}$ at $v$ is parallel to the hyperplane
$\Psi$. The two points $p$ and $v$ of tangency are related by
\begin{equation} \tag{3.31}
  \begin{aligned}
v=\frac{\sqrt{k+h}}{\sqrt{k}}p, \quad p=(p_1, \cdots , p_n,\sqrt{k-r^2}), r^2=a_1^2p_1^2+\cdots+a_n^2p_n^2.
\end{aligned}
  \end{equation}

Then the linear mapping
\begin{equation}\tag{3.32}
   \begin{aligned}
   T_1(x_1,x_2,\cdots , x_n,z)=(a_1x_1,a_2x_2,\cdots , a_nx_n, z)
    \end{aligned}
   \end{equation}
transforms $M_k$ (resp., $M_{k+h}$) onto a hypersphere $M'_k:x_1^2+x_2^2+\cdots +x_n^2+z^2=k,$
(resp., $M_{k+h}':x_1^2+x_2^2+\cdots +x_n^2+z^2=k+h $), $\Phi$ to a hyperplane $\Phi'$,
$p\in M_{k}$ and $v\in M_{k+h}$ to points of tangency $p' =(p_1', \cdots, p_n',\sqrt{k-(r')^2})\in M'_{k}$
and $v' =(\sqrt{k+h}/\sqrt{k})p'\in M'_{k+h}$, respectively,
where $(r')^2=\Sigma (p_i')^2$. The corresponding volume $V_{p'}'^*(k,h)$ is given by
\begin{equation}\tag{3.33}
   \begin{aligned}
   V_{p'}'^*(k,h)=a_1\cdots a_n  V_p^*(k,h).
    \end{aligned}
   \end{equation}
By the symmetry of hyperspheres $M'_k$ and $M'_{k+h}$ centered at the origin, we see that $V_{p'}'^*(k,h)$   is
independent of the point $p'$. This, together with (3.33),  shows that the function $g$ satisfies Condition $(V^*)$.
  \vskip 0.30cm
Finally, we  show that the function $g$ given by $g(x,z)=z^2+a_1^2x_1^2+\cdots +a_n^2x_n^2$
 with  $a_1, \cdots, a_n >0$
satisfies Condition $(A^*)$.

For a fixed point $p=(p_1, \cdots, p_n, \sqrt {k-r^2})\in M_k$,
where $r^2=a_1^2p_1^2+ \cdots +a_n^2p_n^2$, and a  small  $t>0$,
we have $V_p(t)=\phi_k(h(t))$ for some negative function $h=h(t)$ with $h(0)=0$.
 By differentiating with respect to $t$,  we get from (2.5)
 \begin{equation} \tag{3.34}
  \begin{aligned}
  A_p(t)=V_p'(t)=\phi_k'(h)h'(t),
  \end{aligned}
  \end{equation}
where $\phi_k'(h)$ denotes the derivative of $\phi_k$ with respect to $h$.

With the help of (3.31), it is straightforward to show that
the distance $t$ from $p\in M_k$ to the tangent hyperplane $\Phi$ to $M_{k+h}$ at $v\in M_{k+h}$
is given by
 \begin{equation} \tag{3.35}
  \begin{aligned}
 t=\frac{2\sqrt{k}(\sqrt{k}-\sqrt{k+h})}{|\nabla g(p)|}.
  \end{aligned}
  \end{equation}
Hence we get
\begin{equation} \tag{3.36}
  \begin{aligned}
 h(t)=\frac{1}{4k}|\nabla g(p)|^2t^2 - |\nabla g(p)|t.
  \end{aligned}
  \end{equation}
It follows from (3.34), (3.35) and (3.36) that
\begin{equation} \tag{3.37}
  \begin{aligned}
  A_p^*(k,h)=A_p(t)=-\frac{\sqrt{k+h}}{\sqrt{k}}\phi_k'(h)|\nabla g(p)|.
  \end{aligned}
  \end{equation}
This shows that the function $g$ given by $g(x,z)=z^2+a_1^2x_1^2+\cdots +a_n^2x_n^2$
 with  $a_1, \cdots, a_n >0$ satisfies Condition $(A^*)$.

 It follows from (3.28) and (3.29) with $\alpha=2$  that the family $M_k$ of ellipsoids satisfies
   \begin{equation}\tag{3.38}
   \begin{aligned}
 K(p)|\nabla g(p)|^{n+2}=c(k)=2^{n+2}a_1^2a_2^2\cdots a_n^2k.
    \end{aligned}
   \end{equation}

This completes the proof of Theorem 3.

 \vskip 0.2cm \par

\par \vskip 0.2cm
\noindent {\bf 4. Condition $(S^*)$} \vskip0.2cm
\par

In this section, we prove Theorem 5.

Let $f:{\Bbb R}^{n}\rightarrow {\Bbb R}$ be a nonnegative strictly convex function.
 For a real number $\alpha\in R$ with $\alpha\ne 0,1$, we consider
  the function $g:{\Bbb R}^{n+1}\rightarrow {\Bbb R}$  defined by  $g(x,z)=z^\alpha-f(x)$.

 Suppose that  the level hypersurfaces $M_k, k\in S_g$ of $g$  in the $(n+1)$-dimensional Euclidean space
 ${\Bbb E}^{n+1}$ are strictly convex and that the function $g$ satisfies Condition $(S^*)$.
  Then, as in the proof of Theorem 2, we can show that the function $g$ is given by
$g(x,z)=z^2-f(x), z>0$ with $k>0$, where  $f(x)=a_1^2x_1^2+\cdots +a_n^2x_n^2, a_1, \cdots, a_n >0$.

For a fixed $k>0$ and a small $h>0$, consider
 the tangent hyperplane $\Phi$ of $M_{k+h}$ at a point
 $v\in M_{k+h}$ which is  parallel to the tangent hyperplane $\Psi$ of $M_{k}$ at $p\in M_{k}$.
  The two points $p$ and $v$ of tangency are related by
\begin{equation} \tag{4.1}
  \begin{aligned}
v=\frac{\sqrt{k+h}}{\sqrt{k}}p=\frac{\sqrt{k+h}}{\sqrt{k}}(p_1, \cdots , p_n,\sqrt{r^2+k}), r^2=a_1^2p_1^2+\cdots+a_n^2p_n^2.
\end{aligned}
  \end{equation}
The tangent hyperplane $\Phi$ of $M_{k+h}$ at
 $v\in M_{k+h}$ is given by
 \begin{equation} \tag{4.2}
  \begin{aligned}
\Phi: z=\frac{1}{\sqrt{r^2+k}}\{a_1^2p_1x_1+\cdots+a_n^2p_nx_n+\sqrt{k(k+h)}\}.
\end{aligned}
  \end{equation}

The linear transformation $y_1=a_1x_1, \cdots, y_n=a_nx_n, z=z$ transforms $M_k$ to $M_k': z^2=|y|^2+k$,
 $p=(p_1, \cdots , p_n,\sqrt{r^2+k})$ to $q=(a_1p_1, \cdots, a_np_n, \sqrt{r^2+k})$, and  $\Phi$ to the hyperplane
 $\Phi'$ defined by
\begin{equation} \tag{4.3}
  \begin{aligned}
\Phi': z=\frac{1}{\sqrt{|q|^2+k}}\{\<q,y\>+\sqrt{k(k+h)}\}.
\end{aligned}
  \end{equation}
Hence the $n$-dimensional surface area $S_p^*(k,h)$ of the region of $M_k$
between the two hyperplanes $\Phi$ and $\Psi$ is given by
\begin{equation} \tag{4.4}
  \begin{aligned}
S_p^*(k,h)
=\frac{1}{a}\int _{D_q(k,h)} \frac{\{(a_1^2+1)y_1^2+\cdots+(a_n^2+1)y_n^2+k\}^{1/2}}{\sqrt{|y|^2+k}}dy,
\end{aligned}
  \end{equation}
where $a=a_1\cdots a_n$ and
\begin{equation} \tag{4.5}
  \begin{aligned}
D_q(k,h): (|q|^2+k)(|y|^2+k)\le (\<q,y\>+\sqrt{k(k+h)})^2.
\end{aligned}
  \end{equation}

By assumption, $S_p^*(k,h)/|\nabla g(p)|=\eta_k(h)$ is  independent of $p$.
Since we have
$$|\nabla g(p)|^2=4((a_1^2+1)q_1^2+\cdots+(a_n^2+1)q_n^2+k),$$
we see that
\begin{equation} \tag{4.6}
  \begin{aligned}
\int _{D_q(k,h)}H(y)dy=2a\sqrt{|q|^2+k}H(q)\eta_k(h),
\end{aligned}
  \end{equation}
where we denote
\begin{equation} \tag{4.7}
  \begin{aligned}
H(y)= \frac{\{(a_1^2+1)y_1^2+\cdots+(a_n^2+1)y_n^2+k\}^{1/2}}{\sqrt{|y|^2+k}}.
\end{aligned}
  \end{equation}

It is straightforward to show that $D_q(k,h)$ is an ellipsoid centered at $\sqrt{(k+h)/k}q$
and its canonical form is given by
\begin{equation} \tag{4.8}
  \begin{aligned}
\frac{ky_1^2}{h(|q|^2+k)}+\frac{y_2^2+\cdots+y_n^2}{h}\le 1.
\end{aligned}
  \end{equation}
This shows that the volume of $D_q(k,h)$ is given by
\begin{equation} \tag{4.9}
  \begin{aligned}
V( D_q(k,h))=\frac{1}{\sqrt{k}}(\sqrt{h})^n\sqrt{|q|^2+k} \omega_n.
\end{aligned}
  \end{equation}

Let's denote by $\theta_k(h)$  the function defined by
\begin{equation} \tag{4.10}
  \begin{aligned}
\theta_k(h)=\frac{2a\sqrt{k}}{\omega_n(\sqrt{h})^n}\eta_k(h).
\end{aligned}
  \end{equation}
Then, it follows from (4.6) and (4.9) that $H(y)$ satisfies
\begin{equation} \tag{4.11}
  \begin{aligned}
\frac{1}{V( D_q(k,h))}\int _{D_q(k,h)}H(y)dy=H(q)\theta_k(h).
\end{aligned}
  \end{equation}

When $q=0$, $H(0)=1$ and $D_0(k,h)$ is the  ball $B_0(\sqrt{h})$ of radius $\sqrt{h} $ centered at $y=0$.
Hence, (4.11) implies that for any positive numbers $k$ and $h$
\begin{equation} \tag{4.12}
  \begin{aligned}
\theta_k(h)=\frac{1}{V(B_0(\sqrt{h}))}\int _{B_0(\sqrt{h})}H(y)dy >1.
\end{aligned}
  \end{equation}
If we let $a_1=\max \{a_i|i=1, 2, \cdots, n\}$, then we have from (4.7)
\begin{equation} \tag{4.13}
  \begin{aligned}
&1=H(0) \le H(y) < \sqrt{a_1^2+1},\\
&\lim_{y_1\rightarrow\infty}H(y_1, 0, \cdots, 0)=\sqrt{a_1^2+1}.
\end{aligned}
  \end{equation}

Thus,  the left hand side of (4.11) is less than $\sqrt{a_1^2+1}$
for any positive numbers $k, h$ and $q\in {\Bbb R}^{n}$.
But, since $\theta_k(h)>1$, for  $q=(q_1, 0, \cdots, 0)$ with  sufficiently large $q_1$,
the right hand side of (4.11) is greater than
 $\sqrt{a_1^2+1}$. This contradiction completes the proof of Theorem 5.

\vskip 0.2cm \par

\par \vskip 0.2cm
\noindent {\bf 5. Elliptic paraboloids} \vskip0.2cm
\par

In this section, we prove Theorem 6.

Let $f:{\Bbb R}^{n}\rightarrow {\Bbb R}$ be a nonnegative  strictly convex function.
 For a real number $\alpha\in R$ with $\alpha\ne 0,2$, let's consider   the function $g$ defined by  $g(x,z)=z^\alpha-f(x)$.
We suppose that  the level hypersurfaces $M_k, k\in S_g$ of $g$  in the $(n+1)$-dimensional Euclidean space
 ${\Bbb E}^{n+1}$ are strictly convex.

 Suppose that  the function $g$  satisfies the condition 1) or 2) in Theorem 6. Then,
 it follows from Lemma 8 or Lemma 9 that  $g$ satisfies the condition 3) in Theorem 6.

 First, we show that the condition 3) in Theorem 6 implies 4) as follows.
  As in the proof of Theorem 2 in Section 3,
  we can show that if $\alpha\ne 0,1,2$, then (3.12) leads to a contradiction.
  Since $\alpha\ne 0,2$, the remaining case is for $\alpha=1$. In this case, it follows from (3.4) that
   \begin{equation} \tag{5.1}
  \begin{aligned}
\det (f_{ij}(x))=c(k).
\end{aligned}
  \end{equation}
 Hence we see that $\det (f_{ij})$ is a positive constant $c$ with  $c(k)=c$.  Thus, $f(x)$ is a  quadratic polynomial
given  by ([1], [4])
\begin{equation} \tag{5.2}
  \begin{aligned}
f(x_1, \cdots, x_n)=a_1^2x_1^2+\cdots +a_n^2x_n^2, \quad a_1, \cdots, a_n >0.
\end{aligned}
  \end{equation}

  Conversely, suppose that the level hypersurfaces are given by
$M_k:z=f(x)+k, z>0$ with $k>0$, where  $f(x)=a_1^2x_1^2+\cdots +a_n^2x_n^2, a_1, \cdots, a_n >0$.
In this case, we have $R_g=S_g=R$ and $I_k=(k, \infty)$.
From the proof of Theorem 5 in [3], we get
 \begin{equation}\tag{5.3}
   \begin{aligned}
   V_p^*(k,h)=\gamma_nh^{(n+2)/2} \quad \text {and} \quad A_p^*(k,h)=\frac{n+2}{2}\gamma_n|\nabla g(p)|h^{n/2},
    \end{aligned}
   \end{equation}
 where
  \begin{equation}\tag{5.4}
   \begin{aligned}
  \gamma_n=\frac{2\sigma_{n-1}}{n(n+2)a_1a_2\cdots a_n}
    \end{aligned}
   \end{equation}
  and $\sigma_{n-1}$ denotes the surface area of the $(n-1)$-dimensional unit sphere.
  It also follows from (3.3) with $\alpha=1$ that
   \begin{equation}\tag{5.5}
   \begin{aligned}
 K(p)|\nabla g(p)|^{n+2}=2^na_1^2a_2^2\cdots a_n^2.
    \end{aligned}
   \end{equation}

  This completes the proof of Theorem 6.
\vskip 0.50cm

 \noindent {\bf Corollary 10.}
  Let $f:{\Bbb R}^{n}\rightarrow {\Bbb R}$ be a nonnegative  strictly convex function.
 For a nonzero real number $\alpha\in R$, let's denote by  $g$ the function defined by  $g(x,z)=z^\alpha-f(x)$.
 Suppose that $R_g=R$ and  the level hypersurfaces $M_k$ of $g$  in the $(n+1)$-dimensional Euclidean space
 ${\Bbb E}^{n+1}$ are strictly convex for all $k\in R$.  Then the following are equivalent.
\vskip0.3cm

 \noindent 1) The function  $g$ satisfies  Condition $(V^*)$.

\noindent 2) The function  $g$ satisfies Condition $(A^*)$.

\noindent 3)  $K(p)|\nabla g(p)|^{n+2}=c(k)$ is constant on each $M_k$.

\noindent 4) For some positive constants $a_1, \cdots, a_n$,
$$g(x,z)=z-(a_1^2x_1^2+\cdots + a_n^2x_n^2).$$

\noindent {\bf Proof.} Suppose that  the function  $g$ satisfies one of the conditions 1), 2) and 3).
Then as above,  we have  $\alpha=1$ or $\alpha=2$.
In case  $\alpha=2$, $c(k)$ is a nonconstant linear function in $k$.
This contradicts to the positivity of $c(k)$. Hence we have  $\alpha=1$.
Thus,  Theorem 6 completes the proof.  $\square$

\vskip 0.50cm



\end{document}